\documentclass[11pt,letterpaper]{article}

\usepackage{fullpage}
\usepackage{amsmath,amsthm,amssymb}

\long\def\symbolfootnote[#1]#2{\begingroup
\def\thefootnote{\fnsymbol{footnote}}\footnote[#1]{#2}\endgroup}

\begin{document}

\title{Comment on revised version of ``The Hadamard circulant conjecture''}
\author{Robert Craigen \and Jonathan Jedwab}
\date{15 November 2011}
\maketitle

\symbolfootnote[0]{
R.~Craigen is with Department of Mathematics, University of Manitoba, 
Winnipeg MB R3T 2N2, Canada.
\par
J.~Jedwab is with Department of Mathematics, 
Simon Fraser University, 8888 University Drive, Burnaby BC V5A 1S6, Canada.
\par
Both authors are supported by NSERC grants.
\par
Email: 
{\tt craigenr@cc.umanitoba.ca}, {\tt jed@sfu.ca}
\par
2010 Mathematics Subject Classification 15B34 (primary), 94A55 (secondary)
}

\begin{abstract}
The revised version of the claim by Hurley, Hurley and Hurley to have 
proved the circulant Hadamard matrix conjecture is mistaken.
\end{abstract}

\maketitle

In January 2011, Hurley, Hurley and Hurley \cite{hurley-hadamard-lms} claimed to
have proved the circulant Hadamard matrix conjecture, but the proof
was mistaken \cite{craigen-jedwab-circulant-hadamard1}.
In September 2011, a revised version \cite{hurley-hadamard-arxiv} 
of the paper \cite{hurley-hadamard-lms} was posted to the arXiv, with 
the comment that
``This is post publication revision of on-line Bull.\ London Math.\ Soc.\ 
version which changes subsection 3.3.'' We show that the revised version
is also mistaken, by summarising part of the argument 
of \cite{hurley-hadamard-arxiv} and then presenting a counterexample.

A \emph{2-block} is a matrix of the form 
$D = \begin{bmatrix} i & j \\ j & i \end{bmatrix}$ for $i, j \in \{1, -1\}$,
and is \emph{even} if $i=j$ and \emph{odd} if $i = -j$.
Suppose there exists a circulant Hadamard matrix~$H$ of order $4n$. Reorder the
rows and columns of~$H$ to form a $2n \times 2n$ matrix $M$ whose entries
are 2-blocks, as in \cite[p.7]{hurley-hadamard-arxiv},
and write the first row of $M$ as
$\begin{bmatrix} M_0 & M_1 & \dots & M_{2n-1} \end{bmatrix}$.
Then exactly $n$ of the 2-blocks $M_i$ are even, and
\begin{equation}
\label{eqn:mi}
\sum_{\mbox{$i$ : $M_i$ and $M_{i+u}$ are even}} M_i M_{i+u} = 
\begin{bmatrix} 0 & 0 \\ 0 & 0 \end{bmatrix} \quad \mbox{for each $u \ne 0$},
\end{equation}
where all matrix subscripts are reduced modulo~$2n$.
Fix $u \ne 0$. Then from~\eqref{eqn:mi}, for each $i$ such that 
$M_i$ and $M_{i+u}$ are even, we can assign a unique $\ell$ such that 
$M_\ell$ and $M_{\ell+u}$ are even and such that 
$M_\ell M_{\ell+u} = -M_i M_{i+u}$. 
We then also assign $i$ to $\ell$, write
$(i,i+u) \sim (\ell,\ell+u)$, and call the index pairs
$(i, i+u)$ and $(\ell, \ell+u)$ \emph{matching}.

An even 2-block $M_i$ is \emph{symmetric} when the 2-block 
$M_{i+n}$ is also even.  The following argument is given 
\cite[p.8]{hurley-hadamard-arxiv} to claim that
``every even block is symmetric'' when $n > 1$. Suppose, for a contradiction,
that $M_i$ is an even block that is not symmetric. Since $n>1$,
there is an even 2-block $M_{i+u}$ for some $u \ne 0$, and
there must be a pair matching $(i, i+u)$.
In each of five exhaustive cases, this forces the existence 
of a further pair of even 2-blocks $(M_j,M_{j+v})$ for some $j$ and $v$,
where $M_j$ is not symmetric, and there must be a pair matching $(j, j+v)$.
Repeat this procedure. Since this procedure ``cannot continue indefinitely,'' 
we obtain a contradiction.

The following is a counterexample to this claimed procedure, using $n=3$
and only the first of the five specified cases:
\[
(M_0,M_1,M_2,M_3,M_4,M_5) = \left(
\begin{bmatrix} + & + \\ + & + \end{bmatrix},
\begin{bmatrix} + & - \\ - & + \end{bmatrix}, 
\begin{bmatrix} - & - \\ - & - \end{bmatrix}, 
\begin{bmatrix} + & - \\ - & + \end{bmatrix}, 
\begin{bmatrix} - & - \\ - & - \end{bmatrix}, 
\begin{bmatrix} + & - \\ - & + \end{bmatrix} 
\right)
\]
(writing $+$ for $1$ and $-$ for $-1$).
The even 2-blocks are $M_0$, $M_2$, and $M_4$, none of which is symmetric. 
Assign the matchings $(0,2) \sim (2,4)$ and $(0,4) \sim (4,2)$.
Let $i=0$ and $j=2$, and follow the procedure 
of \cite[p.8]{hurley-hadamard-arxiv}. 
Since $(0,2) \sim (2,4)$, there must be a pair matching $(0,4)$.
Then, since $(0,4) \sim (4,2)$, there must be a pair matching $(0,2)$.
However $(0,2)$ already has a matching pair $(2, 4)$, 
so the claimed contradiction does not arise.


\begin{thebibliography}{1}

\bibitem{craigen-jedwab-circulant-hadamard1}
R.\ Craigen and J.\ Jedwab.
\newblock Comment on ``{T}he {H}adamard circulant conjecture''.
\newblock {\tt arXiv:1111.3437v1 [math.CO]}.

\bibitem{hurley-hadamard-lms}
B.\ Hurley, P.\ Hurley, and T.\ Hurley.
\newblock The {H}adamard circulant conjecture.
\newblock {\em Bull.\ London Math.\ Soc.}, 2011.
\newblock doi:10.1112/blms/bdq112.

\bibitem{hurley-hadamard-arxiv}
B.\ Hurley, P.\ Hurley, and T.\ Hurley.
\newblock The {H}adamard circulant conjecture.
\newblock {\tt arXiv:1109.0748v1 [math.RA]}.

\end{thebibliography}
\end{document}